\documentclass{article}%
\usepackage{amsmath}
\usepackage{amsfonts}
\usepackage{amssymb}
\usepackage{graphicx}%
\providecommand{\U}[1]{\protect\rule{.1in}{.1in}}

\newenvironment{proof}[1][Proof]{\noindent\textbf{#1.} }{\ \rule{0.5em}{0.5em}}
\begin{document}

\title{HARMONICITY MODULUS AND APPLICATIONS TO THE APPROXIMATION BY POLYHARMONIC FUNCTIONS}
\author{O. I. KOUNCHEV\\Institute of Mathematics, Bulgarian Academy of Sciences, \\Acad. G. Bonchev St. 8, 1113 Sofia, Bulgaria}
\date{}
\maketitle

\textbf{Abstract.} In the present paper we introduce the notion of harmonicity
modulus and harmonicity $K$-functional and apply these notions to prove a
Jackson type theorem for approximation of continuous functions by polyharmonic
functions. For corresponding results on approximation by polynomials see [3, 7].

\textbf{Key words:} Harmonicity modulus, $K$-functional, polyharmonic
functions, Jackson type theorem.

\bigskip

\begin{center}
\textbf{0. Notions and Notations}

\bigskip
\end{center}

Suppose that $D\subset\mathbb{R}^{n}$ is an open, connected and bounded set
($n\geq2$). We shall work with functions $f$ in the space $HC^{r}\left(
D\right)  ,$ $r\geq0$, consisting of all functions $f$ such that $\Delta
^{r}f,$ the $r^{\text{th }}$ power of the Laplacian, exists and is continuous
in $\overline{D}.$ In the space $HC^{0}\left(  \overline{D}\right)  =C\left(
\overline{D}\right)  $ of functions which are continuous in $\overline{D}$ the
usual norm is%
\[
\left\Vert f\right\Vert :=\max_{x\in\overline{D}}\left\vert f\left(  x\right)
\right\vert .
\]
By $B\left(  x;t\right)  $ we will denote an open ball in $\mathbb{R}^{n}$ :%
\[
B\left(  x;t\right)  :=\left\{  y\in\mathbb{R}^{n}:\left\vert x-y\right\vert
<t\right\}  .
\]
For the function $f$, any point $x\in D,$ and a sufficiently small positive
number $h$ we will consider the spherical mean%
\begin{equation}
\mu_{0}\left(  x,h\right)  :=\mu_{0}\left(  f;x,h\right)  :=\frac{1}%
{\omega_{n}}%
{\displaystyle\int_{\Omega_{\xi}}}
f\left(  x+h\xi\right)  d\omega_{\xi}; \label{1}%
\end{equation}
here $\Omega_{\xi}$ denotes the unit sphere in $\mathbb{R}^{n},$ $\omega_{n}$
denotes its area and $d\omega_{\xi}$ is the area element on $\Omega_{\xi}.$

Further we define the quantity%
\begin{equation}
\Delta_{h}\left(  f;x\right)  :=\mu_{0}\left(  f;x,h\right)  -f\left(
x\right)  . \label{2}%
\end{equation}
Throughout the paper we shall use the symbol $C$ as an universal constant.

\bigskip

\begin{center}
\textbf{1. Harmonicity Modulus}

\bigskip
\end{center}

\textbf{DEFINTION 1.} \emph{The harmonicity modulus of the function f in the
domain }$D$\emph{ is defined by}%
\begin{equation}
\omega^{h}\left(  u\right)  :=\omega^{h}\left(  f;u\right)  :=\sup\left\vert
\Delta_{t}\left(  f;x\right)  \right\vert ,\label{3}%
\end{equation}
\emph{where the }$\sup$\emph{ is taken over }$0<t\leq u,$\emph{ and }$B\left(
x;t\right)  \subset D.$

\bigskip

\textbf{REMARK.} \emph{It is clear that }%
\[
\omega^{h}\left(  f;u\right)  \leq\omega_{1}\left(  f;u\right)
\]
\emph{where }$\omega_{1}$\emph{ is the usual first modulus of continuity (see
[5, 9]).}

It is easy to see that we have the representation
\[
\Delta_{h}\left(  f;x\right)  =\frac{1}{2\omega_{n}}%
{\displaystyle\int_{\Omega}}
\left(  f\left(  x+h\xi\right)  -2f\left(  x\right)  +f\left(  x-h\xi\right)
\right)  d\omega_{\xi}.
\]
This implies
\[
\omega^{h}\left(  f;u\right)  \leq\omega_{2}\left(  f;u\right)
\]
for the usual second modulus of continuity (cf. [5, 9]).

\bigskip

\textbf{PROPOSITION 1.} \emph{For every function }$f$\emph{, continuous in
}$\overline{D}$\emph{, the harmonicity modulus has the following properties:}

\emph{1. }$\lim_{t\longrightarrow0}\omega^{h}\left(  f;t\right)  =0;$

\emph{2. }$\omega^{h}\left(  f;u\right)  $\emph{ is a monotone increasing
funciton;}

\emph{3. for every positive }$u$\emph{ the inequality }%
\[
\omega^{h}\left(  f+g;u\right)  \leq\omega^{h}\left(  f;u\right)  +\omega
^{h}\left(  g;u\right)
\]
\emph{holds;}

\emph{4. for every positive number }$u$\emph{ the inequality }%
\[
\omega^{h}\left(  f;u\right)  \leq2\left\Vert f\right\Vert
\]
\emph{holds.}%

\proof
Property 1) follows from the definition of $\mu_{0}$  and the continuity of
the function $f.$ Properties 2) and 3) are evident. Property 4) follows from
the easy-to-check representation%
\begin{equation}
\Delta_{h}\left(  f;x\right)  =\frac{1}{\omega_{n}}%
{\displaystyle\int_{\Omega}}
\left(  f\left(  x+h\xi\right)  -f\left(  x\right)  \right)  d\omega_{\xi
}.\label{4}%
\end{equation}
%

\endproof

Let us introduce the integral operator $J_{0}$ by%
\begin{equation}
J_{0}\left[  \phi;R\right]  :=%
{\displaystyle\int_{0}^{R}}
\left(  r-r^{n-1}R^{-n+2}\right)  \phi\left(  r\right)  dr \label{5}%
\end{equation}
for $n\geq3,$ and by
\begin{equation}
J_{0}\left[  \phi;R\right]  :=%
{\displaystyle\int_{0}^{R}}
r\log\left(  \frac{R}{r}\right)  \phi\left(  r\right)  dr \label{6}%
\end{equation}
for $n=2.$

Further we will need different forms of the classical Pizzetti formula for the
representation of the spherical means (see [2, 8]).

\bigskip

\textbf{THEOREM 1.} \emph{Let the function }$f$\emph{ have a continuous
Laplacian }$\Delta f$\emph{ in the domain }$D.$\emph{ Then the following
representation holds:}%
\begin{equation}
\mu_{0}\left(  f;x,R\right)  =f\left(  x\right)  +l_{n}J_{0}\left[  \mu
_{0}\left(  \Delta f;x,\cdot\right)  ;R\right]  , \label{7}%
\end{equation}
\emph{where }$l_{n}=\frac{1}{n-2}$\emph{ for }$n\geq3$\emph{ and }$l_{2}=1.$

The remainder can also be written as%
\[
J_{0}\left[  \mu_{0}\left(  \Delta f;x,\cdot\right)  ;R\right]  =\mu
_{0}\left(  \Delta f;x,\vartheta R\right)  J_{0}\left[  1,R\right]
\]
with some number $\vartheta=\vartheta\left(  x;R\right)  $ such that
$0<\vartheta<1.$ Since $J_{0}\left[  1;R\right]  =c_{n}R^{2},$ where
$c_{2}=1/4,$ $c_{n}=\frac{n-2}{2n}$ for $n\geq3,$ we have the representation
\begin{equation}
\mu_{0}\left(  f;x,R\right)  =f\left(  x\right)  +d_{n}R^{2}\Delta f\left(
\xi\right)  ,\label{8}%
\end{equation}
where the opint $\xi=\xi\left(  x,R\right)  \in B\left(  x;R\right)  $ and
$d_{n}=c_{n}l_{n}=\frac{1}{2n}.$

\bigskip

\textbf{PROPOSITION 2.} \emph{Suppose that the function }$f$\emph{, defined
and continuous in }$\overline{D},$\emph{ has a continuous Laplacian }$\Delta
f$\emph{ in }$\overline{D}$\emph{ satisfying the inequality }%
\[
\left\vert \Delta f\left(  x\right)  \right\vert \leq M,\qquad x\in
\overline{D}.
\]
\emph{Then the following inequality holds for every positive number }%
$u$\emph{: }%
\[
\omega^{h}\left(  f;u\right)  \leq Md_{n}u^{2}.
\]

The \emph{proof} follows immediately from Pizzetti's formula (8).

The harmonicity modulus plays a role similar to that of the second modulus of
continuity in the one-dimensional case (see [3]). This is well seen from the
following classical result (cf. [8]).

\bigskip

\textbf{THEOREM 2.} \emph{Let }$u$\emph{ be a function defined and integrable
in the domain }$D$\emph{ in }$R^{n}.$\emph{ Then, if }$t>0,$\emph{ we have}%
\[
\omega^{h}\left(  f;t\right)  =0
\]
\emph{if and only if }$f$\emph{ is harmonic in }$D,$\emph{ i.e. }%
\[
\Delta f\left(  x\right)  =0,\qquad x\in D.
\]

Theorem 2 is the motivation for calling $\omega^{h}$ the harmonicity modulus.
We also recall that harmonic functions are considered to be a multivariate
analogue to the linear functions in one dimension.

\bigskip

\begin{center}
\textbf{2. Harmonicity }$K$\textbf{-Functional}

\bigskip
\end{center}

Here we introduce the notion of harmonicity $K$-functional which provides a
basic tool for studying the important properties of the harmonicity modulus.

\bigskip

\textbf{DEFINITION 2.} \emph{For every function }$f\in C\left(  \overline
{D}\right)  $\emph{ and every number }$t>0$\emph{ we define the harmonicity
}$K$\emph{-functional by}%
\begin{equation}
K^{h}\left(  f;t\right)  :=\inf\left\{  \left\Vert f-g\right\Vert
+t^{2}\left\Vert \Delta g\right\Vert \right\}  , \label{9}%
\end{equation}
\emph{where the infimum is taken over all functions }$g\in HC^{1}\left(
\overline{D}\right)  .$\emph{ }

\bigskip

\begin{center}
\textbf{3. Harmonicity Modulus and Harmonicity K---Functional}

\bigskip
\end{center}

The main technical result of the paper is proved in the present and the next
sections. Roughly speaking, it states that the harmonicity modulus and the
harmonicity $K$-functional are equivalent on compact subdomains of D.

The simple part of the equivalence is the following

\bigskip

\textbf{LEMMA 1.} \emph{Let }$D$\emph{ be an open set in }$\mathbb{R}^{n}%
.$\emph{ For all }$t$\emph{ with }$0<t<\infty$\emph{ and }$f\in C\left(
\overline{D}\right)  ,$\emph{ the inequality}%
\begin{equation}
\omega^{h}\left(  f;t\right)  \leq CK^{h}\left(  f;t\right)  \label{10}%
\end{equation}
\emph{holds with some constant }$C>0.$%

\proof
The proof is based on a standard argument. We split $f=f-g+g,$ and apply
Propositions 1 and 2 to obtain the inequality
\begin{align*}
\omega^{h}\left(  f;t\right)   &  \leq\omega^{h}\left(  f-g;t\right)
+\omega^{h}\left(  g;t\right)  \leq2\left\Vert f-g\right\Vert +d_{n}%
t^{2}\left\Vert \Delta g\right\Vert \\
&  \leq\max\left\{  2,d_{n}\right\}  \left[  \left\Vert f-g\right\Vert
+t^{2}\left\Vert \Delta g\right\Vert \right]  .
\end{align*}

Since $g\in HC^{1}\left(  \overline{D}\right)  $ is arbitrary, the statement
(10) is proved. %

\endproof

The domination of $K^{h}$ by $\omega^{h}$ will be established only on compact
subdomains of $D$ in the sense that the $K$-functional of the subdomain
$D_{1},$ $K_{d_{1}}^{h}$ will be proved to be dominated by the harmonicity
modulus $\omega_{D}^{h}$ with respect to the domain $D.$

The problem is that for every $R$ (possibly such that $R<R_{1}$ for some
sufficiently small positive number $R_{1}$) we have to find a function
$g_{R}\in HC^{1}\left(  \overline{D}\right)  $ such that%
\begin{equation}
\left\Vert f-g_{R}\right\Vert +R^{2}\left\Vert \Delta g_{R}\right\Vert \leq
C\omega^{h}\left(  f;R\right)  , \label{11}%
\end{equation}
where the constant $C$ does not depend on $f$ and $R.$

Following the scheme given in [5], taking some spherical means of the function
$f,$ we succeed in constructing the function $g_{R}$ not on the whole of $D$
but on every subdomain $D_{1}$, such that $\overline{D_{1}}\subset D$ and
$R_{1}\leq dist\left(  D_{1},\partial D\right)  .$ In such a way we can prove
the inequality (11) over subdomains where the norm $\left\Vert \cdot
\right\Vert $ is in fact $\left\Vert \cdot\right\Vert _{D_{1}}.$

\begin{center}
\bigskip

\textbf{4. Domination of }$K^{h}$ \textbf{by }$\omega^{h}$ \textbf{ on Compact
Subdomains}

\bigskip
\end{center}

Having in mind Pizzetti's formula (7) in Theorem 1, we consider the function%
\begin{equation}
g_{R,t}\left(  x\right)  =v\left(  t\right)  J_{0,s}\left[  \mu_{0}\left(
f;x,Rs\right)  ;t\right]  . \label{12}%
\end{equation}
Here $J_{0}$ is the operator given by (5) and (6) and $J_{0,s}$ means that $s$
is the input variable for $J_{0}$; the output variable is $t$; $v(t)$ is equal
to $\left(  J_{0}\left[  1;t\right]  \right)  ^{-1},$ where $J_{0}\left[
1;t\right]  $ is the value of the functional for $\phi\left(  t\right)  =1$,
so in fact $\frac{1}{v\left(  t\right)  }=t^{2}\left(  \frac{1}{2}-\frac{1}%
{n}\right)  $ for $n\geq3,$ $\frac{1}{v\left(  t\right)  }=\frac{t^{2}}{4}$
for $n=2.$

The operator $J_{0}$ changes the output in a specific way described by

\bigskip

\textbf{PROPOSITION 3.} \emph{For every integrable function }$\phi$\emph{ and
positive numbers }$s$\emph{ and }$R$\emph{ we have}%
\begin{equation}
J_{0,s}\left[  \phi\left(  st\right)  ;R\right]  =\frac{1}{s^{2}}%
J_{0,t}\left[  \phi\left(  t\right)  ;sR\right]  . \label{13}%
\end{equation}
%

\proof
We give the proof for $n\geq3$. Then $J_{0}$ is given by formula (5). The case
$n=2$ is similar.

By changing the variables we obtain%
\begin{align*}
J_{0,s}\left[  \phi\left(  st\right)  ;R\right]   &  =%
{\displaystyle\int_{0}^{R}}
\left(  t-t^{n-1}R^{-n+2}\right)  \phi\left(  st\right)  dt\\
&  =\frac{1}{s^{2}}%
{\displaystyle\int_{0}^{sR}}
\left(  t-t^{n-1}\left(  sR\right)  ^{-n+2}\right)  \phi\left(  t\right)  dt\\
&  =\frac{1}{s^{2}}J_{0,t}\left[  \phi\left(  t\right)  ;sR\right]  .
\end{align*}%
\endproof

Proposition 3 shows that (12) becomes%
\begin{equation}
g_{R,t}\left(  x\right)  =v\left(  t\right)  \frac{1}{R^{2}}J_{0}\left[
\mu_{0}\left(  f;x,\cdot\right)  ;tR\right]  . \label{14}%
\end{equation}

The following is the main technical result of the paper.

\bigskip

\textbf{THEOREM 3.} \emph{For every subdomain }$D_{1}$\emph{ such that
}$\overline{D_{1}}\subset D,$\emph{ the inequality}%
\begin{equation}
K^{h}\left(  f;R\right)  _{D_{1}}\leq C\omega^{h}\left(  f;R\right)  _{D}
\label{15}%
\end{equation}
\emph{holds for every number }$R$\emph{ with}%
\begin{equation}
0<R<d=dist\left(  D_{1},\partial D\right)  ; \label{16}%
\end{equation}
\emph{here the constant }$C$\emph{ does not depend on }$f$\emph{ and }%
$R$\emph{, and }$K^{h}\left(  f;R\right)  _{D_{1}}$\emph{ denotes the
harmonicity }$K$\emph{-functional for the domain }$D_{1}$\emph{, while
}$\omega^{h}\left(  f;R\right)  _{D}$\emph{ denotes the harmonicity modulus on
}$D.$%

\proof
Let us notice first that the function $g_{R,t}$ given by (12), is well defined
in $D_{1}$ for $R$ satisfying (16) and every number $t$ with $0<t<1.$

Since $v\left(  t\right)  J_{0}\left[  1;t\right]  =1,$ we obtain%
\[
g_{R,t}\left(  x\right)  -f\left(  x\right)  =v\left(  t\right)
J_{0,s}\left[  \mu_{0}\left(  f;x,Rs\right)  -f\left(  x\right)  ;t\right]
\]
for every number $R$ satisfying (16) and every number $t$ with $0\leq t\leq1.$
Hence, for every $x\in D_{1}$ we obtain the inequalities%
\begin{equation}
\left\vert g_{R,t}\left(  x\right)  -f\left(  x\right)  \right\vert \leq
v\left(  t\right)  J_{0}\left[  1;t\right]  \omega^{h}\left(  f;Rt\right)
\leq\omega^{h}\left(  f;R\right)  . \label{17}%
\end{equation}
Consequently, we have proved the domination of the first term of $K^{h}$ by
$\omega^{h}$:
\[
\left\Vert g_{R,t}-f\right\Vert \leq\omega^{h}\left(  f;R\right)  .
\]
For proving the domination of the second term, we will check the value of
$\Delta g_{R,t}\left(  x\right)  $ for $x\in D_{1}.$

First let us suppose that $f$ is twice differentiable in $D,$ i.e. $f\in
C^{2}\left(  D\right)  .$ By formula (14) the Laplacian of $g_{R,t}\left(
x\right)  $ is then equal to
\begin{equation}
\Delta g_{R,t}\left(  x\right)  =v\left(  t\right)  \frac{1}{R^{2}}%
J_{0,s}\left[  \mu_{0}\left(  \Delta f;x,s\right)  ;tR\right]  . \label{18}%
\end{equation}
Hence, combining with formula (7), we obtain
\begin{equation}
\Delta g_{R,t}\left(  x\right)  =v\left(  t\right)  \frac{1}{l_{n}R^{2}%
}\left[  \mu_{0}\left(  f;x,Rt\right)  -f\left(  x\right)  \right]  .
\label{19}%
\end{equation}
This implies the inequalities%
\begin{equation}
R^{2}\left\vert \Delta g_{R,t}\left(  x\right)  \right\vert \leq v\left(
t\right)  \frac{1}{l_{n}}\omega^{h}\left(  f;tR\right)  \leq C\omega
^{h}\left(  f;R\right)  \label{20}%
\end{equation}
for every $x\in D_{1}$ and every $R<dist(D_{1},\partial D),$ where $C$ is a
constant given by
\[
C:=\left\vert v\left(  t_{1}\right)  \right\vert /\not l  _{n}%
\]
and $t_{1}$ is an arbitrary number with $0<t_{1}<1$ satisfying $v\left(
t_{1}\right)  \neq0.$

Inequality (20) implies that the second term of $K^{h}$ is dominated by
$\omega^{h}.$ This ends the proof for $f\in C^{2}\left(  D\right)  .$

In the case of an arbitrary continuous function $f,$ let us take an
approximation to $f$, say $f_{\delta}\in C^{\infty},$ $\delta>0,$ such that
$f_{\delta}$ converges to $f$ uniformly on $D_{1}$ $\delta\longrightarrow0$
(see this construction in [1, paragraph 5]).

Since $\overline{D_{1}}$ is a compact set and $f$ is continuous in $D,$ we
obtain by a standard limiting argument that formula (7) takes for an arbitrary
integrable function $f$ the form%
\[
\mu_{0}\left(  f;x,R\right)  =f\left(  x\right)  +l_{n}\Delta J_{0}\left[
\mu_{0}\left(  f;x,\cdot\right)  ;R\right]  .
\]
This implies that the relations (19) and (20) hold as well. %

\endproof

Now we are ready to prove an important property of the harmonicity modulus.

For an arbitrary subdomain $D_{1}$ of $D$ let us denote by $\omega^{h}\left(
\cdot;\cdot\right)  _{D_{1}}$ the harmonicity modulus for the set $D_{1}.$

\bigskip

\textbf{THEOREM 4.} \emph{For every subdomain }$D_{1}$\emph{ of }$D$\emph{
such that }$\overline{D_{1}}\subset D,$\emph{ the following inequalities hold:
}

\emph{(i) }$\qquad\qquad\qquad\omega^{h}\left(  f;\lambda R\right)  _{D_{1}%
}\leq C\left(  \lambda+1\right)  ^{2}\omega^{h}\left(  f;r\right)  _{D}$\emph{
}

\emph{for every number }$R\leq d=dist(D_{1},\partial D)$\emph{ and every
number }$\lambda>0$\emph{ such that }$\lambda R\leq d;$

\emph{(ii) }$\qquad\qquad\qquad\omega^{h}\left(  f;a+b\right)  _{D_{1}}\leq
C\left[  \omega^{h}\left(  f;a\right)  _{D}+\omega^{h}\left(  f;b\right)
_{D}\right]  $\emph{ }

\emph{for all positive real numbers }$a$\emph{ and }$b.$%

\proof
Inequality (i) follows from a similar inequlity for the harmonicity
$K$-functional. Indeed, since%
\[
\left\Vert f-g\right\Vert _{D_{1}}+\left(  \lambda u\right)  ^{2}\left\Vert
\Delta g\right\Vert _{D_{1}}\leq\left(  \lambda+1\right)  ^{2}\left(
\left\Vert f-g\right\Vert _{D_{1}}+u^{2}\left\Vert \Delta g\right\Vert
_{D_{1}}\right)
\]
for an arbitrary function $g\in HC^{1}\left(  \overline{D_{1}}\right)  ,$ by
the definition of the harmonicity $K$-functional we obtain the inequality%
\[
K^{h}\left(  f;\lambda u\right)  _{D_{1}}\leq\left(  \lambda+1\right)
^{2}K^{h}\left(  f;u\right)  _{D_{1}},
\]
for every number $\lambda\geq0.$

Lemma 1 gives
\[
\omega^{h}\left(  f;\lambda R\right)  _{D_{1}}\leq CK^{h}\left(  f;\lambda
R\right)  _{D_{1}},
\]
and Theorem 3 implies%
\[
K^{h}\left(  f;\lambda u\right)  _{D_{1}}\leq C\omega^{h}\left(  f;R\right)
_{D}%
\]
for every number $R<d$. These inequalities imply the inequality (i).

In order to prove inequality (ii) let us note that%
\[
\left(  a+b\right)  ^{2}\leq2\left(  a^{2}+b^{2}\right)
\]
for all real numbers $a$ and $b.$ This implies%
\begin{align*}
&  \left\Vert f-g\right\Vert _{D_{1}}+\left(  a+b\right)  ^{2}\left\Vert
\Delta g\right\Vert _{D_{1}}\\
&  \leq\left\Vert f-g\right\Vert _{D_{1}}+2\left(  a^{2}+b^{2}\right)
\left\Vert \Delta g\right\Vert _{D_{1}}\\
&  \leq2\left\Vert f-g\right\Vert _{D_{1}}+2a^{2}\left\Vert \Delta
g\right\Vert _{D_{1}}+2\left\Vert f-g\right\Vert _{D_{1}}+2b^{2}\left\Vert
\Delta g\right\Vert _{D_{1}}.
\end{align*}
The definition of $K^{h}$ implies the inequality%
\[
K^{h}\left(  f;a+b\right)  _{D_{1}}\leq2\left(  K^{h}\left(  f;a\right)
_{D_{1}}+K^{h}\left(  f;b\right)  _{D_{1}}\right)  .
\]
Now inequality (ii) follows by arguments similar to those used for (i).%

\endproof

\begin{center}
\bigskip\textbf{5. Polyharmonic Kernels}
\end{center}

Here we introduce kernels which are polyharmonic functions and arise naturally
from the Jackson type kernels used in approximation theory [3].

Let us recall that the function $f$ is called polyharmonic of order $p$ in an
open set $D,$ where $p$ is a nonnegative integer, if it satisfies the equation%
\[
\Delta^{p}f\left(  x\right)  =0\qquad\text{for }x\in D;
\]
here the iterated Laplacian of order $p$ is defined inductively by the
equations $\Delta^{k+1}:=\Delta\Delta^{k}$ for $k\geq0$ and $\Delta^{0}:=id$
(see [8]).

Let us remind the notion of Jackson type kernel (cf. [3]).

\bigskip

\textbf{DEFINITION 3.} \emph{A kernel of Jackson type of order }$\nu$\emph{,
where }$\nu=1,2,...,$\emph{ is defined to be the function given by}%
\[
J_{k;\nu}\left(  t\right)  :=\left(  \gamma_{k,\nu}\right)  ^{-1}\left[
\sin\left(  \nu t/2\right)  /\sin\left(  t/2\right)  \right]  ^{2k},
\]
\emph{where }$k$\emph{ is a natural number and the constant is}%
\[
\gamma_{k,\nu}:=\frac{1}{\pi}%
{\displaystyle\int_{-\pi}^{\pi}}
\left[  \sin\left(  \nu t/2\right)  /\sin\left(  t/2\right)  \right]
^{2k}dt.
\]

For the properties of these kernels we refer to [3]. Through the substitution%
\[
x=2\sin\left(  t/2\right)  ,\qquad t\in\left[  -\pi,\pi\right]  ,\qquad
x\in\left[  -2,2\right]  ,
\]
we obtain the nonperiodic Jackson type kernels:
\[
\overline{J}_{k,\nu}\left(  x\right)  =\gamma_{k,\nu}\left(  \overline{\gamma
}_{k,\nu}\right)  ^{-1}J_{k,\nu}\left[  \arccos\left(  1-x^{2}/2\right)
\right]  ;
\]
here the constant is%
\[
\overline{\gamma}_{k,\nu}:=%
{\displaystyle\int_{-1}^{1}}
\gamma_{k,\nu}J_{k,\nu}\left[  \arccos\left(  1-x^{2}/2\right)  \right]  dx,
\]
for $\nu\in\mathbb{N}.$

Finally, we define the polyharmonic Jackson type kernels of order $p$ by the
equation%
\[
\widetilde{J}_{k,p}\left(  x\right)  :=\widetilde{J}_{k,p}\left(  \left\vert
x\right\vert \right)  =\overline{\gamma}_{k,p}\left(  \widetilde{\gamma}%
_{k,p}\right)  ^{-1}\overline{J}_{k,p}\left(  \left\vert x\right\vert \right)
,
\]
for $p\in\mathbb{N},$ and for every $x\in\mathbb{R}^{n}$ such that $\left\vert
x\right\vert \leq2;$ here the constant is given by
\[
\widetilde{\gamma}_{k,p}=%
{\displaystyle\int_{0}^{1}}
r^{n-1}\overline{\gamma}_{k,p}\overline{J}_{k,p}\left(  r\right)  dr=%
{\displaystyle\int_{0}^{1}}
r^{n-1}\left[  \sin\left(  \nu t/2\right)  /\sin\left(  t/2\right)  \right]
^{2k}dr,
\]
where $t=\arccos\left(  1-r^{2}/2\right)  .$ 

\bigskip

\textbf{THEOREM 5.} \emph{The polyharmonic Jackson type kernels have the
following properties:}

\emph{(i) For all natural numbers }$p$\emph{ and }$k,$\emph{ the kernel
}$\widetilde{J}_{k,p}$\emph{ is a nonnegative polyharmonic function of order
}$k(p-1)+1$\emph{ and }$\widetilde{J}_{k,p}\left(  x\right)  $\emph{ is
defined for every }$x\in\mathbb{R}^{n}$\emph{ satisfying }$\left\vert
x\right\vert \leq2$\emph{;}

\emph{(ii) }$%
{\displaystyle\int_{B\left(  0;1\right)  }}
\widetilde{J}_{k,p}\left(  x\right)  =1$\emph{ ;}

\emph{(iii) If }$I_{i}$\emph{ is defined by }%
\[
I_{i}:=%
{\displaystyle\int_{0}^{1}}
t^{i+n-1}\widetilde{J}_{k,p}\left(  t\right)  dt
\]
\emph{for nonnegative integers }$i,$\emph{ then for }$i<2k-n$\emph{ we have
the inequality }%
\[
I_{i}\leq Cp^{-i},
\]
\emph{and for }$i=2k-n$\emph{ we have the inequality}
\[
I_{i}\leq C\left(  \ln p\right)  p^{-i}.
\]

The proof of Theorem 5 is based on standard arguments [3] and is given in
detail in a forthcoming paper [6].

\bigskip

\begin{center}
\textbf{6. A Direct Theorem of Jackson Type}

\bigskip
\end{center}

Here we prove an approximation theorem which is analogous to the direct
theorem of Jackson for the approximation by polynomials in the one-dimensional
case, where the rate of approximation is estimated by the first and the second
modulus of continuity (see [3, 7]).

In the multivariate case we approximate by polyharmonic functions and the rate
of approximation is estimated by the harmonicity modulus.

Let us first give some necessary notations. Let $K$ be a polyharmonic function
on $\left\{  x\in\mathbb{R}^{n}:\left\vert x\right\vert \leq1\right\}  .$ Then
for every function $f$ defined and continuous in the domain $D,$ we can define
the operator%
\begin{equation}
T_{K}\left[  f\right]  \left(  x\right)  :=%
{\displaystyle\int_{B\left(  x;1\right)  }}
K\left(  x-u\right)  f\left(  u\right)  du \label{21}%
\end{equation}
for every $x\in D$ such that $dist(x,\partial D)<1.$

Let the domain $D$ be regular in the sense of solvability of the Dirichlet
problem (see [4]), and let the function $f$ be continuous in $\overline{D}.$
Then there exists a harmonic function $h_{f}$ solving the Dirichlet problem in
$D,$ i.e.%
\begin{align*}
\Delta h_{f}\left(  x\right)   &  =0\qquad\text{for }x\in D,\\
h_{f}\left(  x\right)   &  =f\left(  x\right)  \qquad\text{for }x\in\partial
D.
\end{align*}
We shall consider the function $F_{0}$ given by the following conditions:%
\[
F_{0}\left(  x\right)  :=f\left(  x\right)  -h_{f}\left(  x\right)
\qquad\text{for }x\in\overline{D},
\]
and
\[
F_{0}\left(  x\right)  :=0\qquad\text{for }x\notin\overline{D}.
\]
The function $F_{0}$ is evidently continuous on the whole space and it makes
sense to consider its harmonicity modulus there or in domains containing $D.$

Another interesting feature of the function $F_{0}$ is that its harmonicity
modulus in $D$ satisfies%
\[
\omega^{h}\left(  F_{0};t\right)  _{D}=\omega^{h}\left(  f;t\right)  _{D}.
\]
This follows immediately from the Gauss mean value theorem, which states that%
\[
\mu_{0}\left(  h_{f};x,t\right)  =h_{f}\left(  x\right)
\]
for every $x\in D$ and $t>0$ such that $B(x;t)\subset D.$

Notice that for every domain $D_{1}$ such that $\overline{D}\subset D_{1}$ we
have%
\[
\omega^{h}\left(  F_{0};t\right)  _{D_{1}}=\omega^{h}\left(  F_{0};t\right)
_{\mathbb{R}^{n}}%
\]
for every positive number $t\leq dist(D,\partial D_{1}).$ Here $\omega
^{h}\left(  F_{0};t\right)  _{\mathbb{R}^{n}}$ denotes the harmonicity modulus
of the function $F_{0}$ in the whole space.

Let $D_{2}$ be a domain such that $\overline{D_{1}}\subset D_{2}.$ Then we can
apply Theorem 3 to obtain the following inequalities%
\begin{equation}
C_{1}\omega^{h}\left(  F_{0};t\right)  _{D}\leq K^{h}\left(  F_{0};t\right)
_{D}\leq C_{2}\omega^{h}\left(  F_{0};t\right)  _{D_{2}} \label{22}%
\end{equation}
for sufficiently small numbers $t>0$ and appropriate constants $C_{1},$
$C_{2}$ which do not depend on $f$ and $t.$

Next suppose that, for some nonnegative $r,$ the function $f$ is in
$HC^{r}(\overline{D}).$ Then, inductively in $r,$ we obtain a solution $h_{f}$
to the following boundary value problem:%
\begin{align*}
\Delta^{r+1}h_{f}\left(  x\right)   &  =0,\qquad\qquad x\in D;\\
\Delta^{j}h_{f}\left(  x\right)   &  =\Delta^{j}f\left(  x\right)  ,\qquad
x\in\partial D,
\end{align*}
for $j=0,1,...,r.$

We shall consider the function $F_{r}$ given by
\begin{align}
F_{r}\left(  x\right)   &  :=f\left(  x\right)  -h_{f}\left(  x\right)
\qquad\text{for }x\in\overline{D};\label{23}\\
F_{r}\left(  x\right)   &  :=0\qquad\qquad\qquad\text{for }x\notin\overline
{D}.\nonumber
\end{align}
Note that the function $F_{r}$ is continuous on the whole space together with
$\Delta^{r}F_{r}$ and we can apply to it all properties of the harmonicity
modulus, a fact which will be used below.

Now we are ready to state the following result which is the main application
of the harmonicity modulus in the present paper.

\bigskip

\textbf{THEOREM 6.} \emph{Let the domain }$D$\emph{ be regular in the sense of
solvability of the Dirichlet problem. Let for some integer }$r\geq0,$\emph{
the function }$f\in HC^{r}\left(  \overline{D}\right)  .$\emph{ Let us denote
by }$F_{r}$\emph{ the function given by (23). Then, for every natural number
}$p$\emph{ satisfying }$p\geq r+1$\emph{, there exists a polyharmonic function
}$T_{p}$\emph{ of order }$p$\emph{ in }$D$\emph{ satisfying the inequality}%
\begin{equation}
\left\vert f\left(  x\right)  -T_{p}\left(  x\right)  \right\vert \leq
C\omega^{h}\left(  \Delta^{r}F_{r};\frac{1}{p}\right)  \frac{1}{p^{2r}}
\label{24}%
\end{equation}
\emph{for every }$x\in\overline{D},$\emph{ where the constant }$C>0$\emph{
depends on the domain }$D$\emph{ and on }$r.$%

\proof
(1) By a similarity transform we can suppose that the domain $D$ is contained
in the ball $B(0;1/2).$ Obviously, this transform preserves the polyharmonic
functions. To find the harmonicity modulus for the function $G(x)=f(\lambda
x),$ where $\lambda$ is a positive real number, let us compute the harmonicity
difference given by (2):%
\[
\Delta_{t}\left(  G;x\right)  =\mu_{0}\left(  G;x,t\right)  -G\left(
x\right)  =\mu_{0}\left(  f;\lambda x;\lambda t\right)  -f\left(  \lambda
x\right)  =\Delta_{\lambda t}\left(  f;\lambda x\right)  .
\]
Hence, we obtain
\[
\omega^{h}\left(  G;t\right)  _{D}=\omega^{h}\left(  f;\lambda t\right)
_{\lambda D},
\]
where $\lambda D$ is the domain given by
\[
\lambda D=\left\{  y\in\mathbb{R}^{n}:y=\lambda x,\quad x\in D\right\}  .
\]
So for a domain $D_{1}$ such that $\overline{D_{1}}\subset D$ we have
\[
\omega^{h}\left(  G;t\right)  _{D_{1}}=\omega^{h}\left(  f;\lambda t\right)
_{\lambda D_{1}},
\]
which proves the inequality
\[
\omega^{h}\left(  G;t\right)  _{D_{1}}\leq\left(  \lambda+1\right)  ^{2}%
\omega^{h}\left(  f;t\right)  _{\lambda D}%
\]
for every number $t\leq dist\left(  \lambda D_{1},\partial\lambda D\right)  .$

This shows that the harmonicity modulus is at most multiplied by a constant as
a result of a similarity transform. Applied to the modulus $\omega^{h}\left(
\Delta^{r}F_{0};p^{-1}\right)  _{\mathbb{R}^{n}}$ we see that by (22) it only
changes up to a constant multiple.

(2)\qquad We will define the polyharmonic function $T_{p}\left(  x\right)
=T_{p}\left(  f;r,x\right)  $ of order $p$ inductively by the following
recurrency relation:
\begin{align}
T_{p}\left(  x\right)   &  :=T_{p}\left(  F_{r};m,x\right)  \label{25}\\
&  :=T_{p}\left(  F_{r};m-1,x\right)  +T_{k,\nu}\left[  F_{r}\left(
\cdot\right)  -T_{p}\left(  F_{r};m-1,\cdot\right)  \right]  \left(  x\right)
,\nonumber
\end{align}
for every $x\in\overline{D}$ and every $m$ with $1\leq m\leq r.$

Here $T_{k,\nu}$ is a short notation for the operator given by formula (21)
for the Jackson Type kernel $\widetilde{J}_{k,\nu},$ where we take $k$ big
enough to satisfy $2k-n\geq3$, and put $\nu:=[(p-1)/k]+1$ (here $[y]$ denotes,
as usually, the greatest integer which does not exceed $y$). The choice of
such $\nu$ provides that the order of the polyharmonic function $\widetilde
{J}_{k,p}$ be equal to $k(\nu-1)+1<p.$

Note that the operator $T_{k,\nu}$ is well defined and produces a polyharmonic
function since $F_{r}$ is a finite function, the kernels are defined in
$B(0;1)$ and we have the inclusion $D\subset B\left(  0;1/2\right)  .$

(3)\qquad Let us check the Theorem for $r=0.$ In this case we have $f\in
C\left(  \overline{D}\right)  .$

Due to Theorem 5 the following holds:%
\begin{align*}
D\left(  x\right)   &  :=R_{0}\left(  x\right)  -T_{k,\nu}\left[
F_{0}\right]  \left(  x\right) \\
&  =%
{\displaystyle\int_{B\left(  x;1\right)  }}
\left[  F_{0}\left(  x\right)  -F_{0}\left(  u\right)  \right]  \widetilde
{J}_{k,\nu}\left(  x-u\right)  du\\
&  =%
{\displaystyle\int_{0}^{1}}
\left\{
{\displaystyle\int_{\Omega_{\xi}}}
\left[  F_{0}\left(  x\right)  -F_{0}\left(  x-r\xi\right)  \right]
d\omega_{\xi}\right\}  r^{n-1}\widetilde{J}_{k,\nu}\left(  r\right)  dr.
\end{align*}
By the properties of the harmonicity modulus (see Theorem 4) this gives the
following estimate%
\begin{align*}
\left\vert D\left(  x\right)  \right\vert  &  \leq\omega_{n}%
{\displaystyle\int_{0}^{1}}
r^{n-1}\widetilde{J}_{k,\nu}\left(  r\right)  \omega^{h}\left(  F_{0}%
;r\right)  _{\mathbb{R}^{n}}dr\\
&  \leq C\omega^{h}\left(  F_{0};p^{-1}\right)  _{\mathbb{R}^{n}}\omega_{n}%
{\displaystyle\int_{0}^{1}}
r^{n-1}\widetilde{J}_{k,\nu}\left(  r\right)  \left(  pr+1\right)  ^{2}dr
\end{align*}
for every $p\geq1$ and some constant $C>0.$

Again, applying Theorem 5, (iii), since $2k-n\geq3$, we have the inequality%
\[%
{\displaystyle\int_{0}^{1}}
r^{n-1}\widetilde{J}_{k,\nu}\left(  r\right)  \left(  pr+1\right)  ^{2}dr\leq
C%
{\displaystyle\int_{0}^{1}}
r^{n-1}\widetilde{J}_{k,\nu}\left(  r\right)  \left(  \nu r+br+1\right)
^{2}dr\leq C_{1}%
\]
for appropriate constants $C,$ $C_{1}$ and $b.$The last gives, finally, that%
\begin{equation}
\left\vert D\left(  x\right)  \right\vert \leq C\omega^{h}\left(  F_{0}%
;p^{-1}\right)  ,\qquad x\in\overline{D}, \label{26}%
\end{equation}
for some constant $C>0.$ From this estimate we get the statement for $r=0.$

(4)\qquad Before proceeding by induction on $r,$ let us note the following. If
for some function $\phi$ on some domain $D,$ such that $\Delta\phi$ is
continuous on $D,$ the inequality%
\[
\left\vert \Delta\phi\left(  x\right)  \right\vert \leq M,\qquad x\in D,
\]
holds, then by Proposition 2 we obtain the inequality%
\[
\omega^{h}\left(  \phi;t\right)  \leq Md_{n}t^{2}%
\]
for every number $t>0.$ Hence, by (26), we obtain the inequality%
\begin{equation}
\left\vert \phi\left(  x\right)  -T_{k,\nu}\left[  \phi\right]  \left(
x\right)  \right\vert \leq CM\frac{1}{p^{2}}\label{27}%
\end{equation}
for an appropriate constant $C.$

(5)\qquad Let us suppose that the Theorem is true for the classes of functions
$HC^{0},$ $HC^{1},$ $...,$ $HC^{r},$ $r\geq0.$ Then, if $f\in HC^{r+1},$ it
follows that $\Delta f\in HC^{r},$ and equality (25) implies that%
\[
\Delta T_{p}\left(  F_{r+1};r,x\right)  =T_{p}\left(  \Delta F_{r+1}%
;r,x\right)  .
\]
Applied to the function $\Delta F_{r+1},$ the induction hypothesis (24) gives
\begin{align*}
\left\vert \Delta\left[  F_{r+1}\left(  x\right)  -T_{p}\left(  F_{r+1}%
;r,x\right)  \right]  \right\vert  &  =\left\vert \Delta F_{r+1}\left(
x\right)  -T_{p}\left(  \Delta F_{r+1};r,x\right)  \right\vert \\
&  \leq C\omega^{h}\left(  \Delta^{r+1}F_{r+1};p^{-1}\right)  _{D}\times
p^{-2r}.
\end{align*}
Let us put
\[
\phi\left(  x\right)  :=F_{r+1}\left(  x\right)  -T_{p}\left(  F_{r+1}%
;r,x\right)
\]
and apply inequality (27) to this function $\phi.$ We obtain the following
inequalities:%
\begin{align}
&  \left\vert \phi\left(  x\right)  -T_{k,\nu}\left[  \phi\right]  \left(
x\right)  \right\vert \label{28}\\
&  =\left\vert F_{r+1}\left(  x\right)  -T_{p}\left(  F_{r+1};r,x\right)
-T_{k,\nu}\left[  F_{r+1}\left(  \xi\right)  -T_{p}\left(  F_{r+1}%
;r;\xi\right)  \right]  \left(  x\right)  \right\vert \nonumber\\
&  \leq CC_{1}p^{-2}\omega^{h}\left(  \Delta^{r+1}F_{r+1};p^{-1}\right)
_{D}\frac{1}{p^{2r}}\nonumber\\
&  =CC_{1}p^{-2}\omega^{h}\left(  \Delta^{r+1}F_{r+1};p^{-1}\right)  _{D}%
\frac{1}{p^{2\left(  r+1\right)  }}.\nonumber
\end{align}
On the other hand, by (25) we have%
\[
\phi\left(  x\right)  -T_{k,\nu}\left(  \phi\right)  \left(  x\right)
=F_{r+1}\left(  x\right)  -T_{p}\left(  F_{r+1};r+1,x\right)  ,
\]
which shows that the inequality in (28) is exactly inequality (24) for $r+1.$
This yields the statement of the Theorem for $r+1.$ %

\endproof

\bigskip

\textbf{COROLLARY.} \emph{In view of the Remark after Definition 1, in Theorem
6 we can replace inequality (24) by the following inequalities}%
\[
\left\vert f\left(  x\right)  -T_{p}\left(  x\right)  \right\vert \leq
C\omega_{1}\left(  \Delta^{r}F_{r};\frac{1}{p}\right)  \frac{1}{p^{2r}}%
\]
\emph{or}%
\[
\left\vert f\left(  x\right)  -T_{p}\left(  x\right)  \right\vert \leq
C\omega_{2}\left(  \Delta^{r}F_{r};\frac{1}{p}\right)  \frac{1}{p^{2r}}%
\]
\emph{for }$x\in\overline{D},$\emph{ where }$\omega_{1}$\emph{ and }%
$\omega_{2}$\emph{ are the usual first and second moduli of continuity (see
[5]).}

\bigskip

\begin{center}
\textbf{Acknowledgements}
\end{center}

The author wishes to thank Professor Walter K. Hayman and the reviewer for
their invaluable remarks.

\begin{center}
\textbf{References}
\end{center}

1.\qquad Besov, O. V., V. P. Hin and S. M. Nikolskii: \emph{Integral
Representation of Functions and Imbedding Theorems}, Nauka, Moskow 1975.

2.\qquad Courant, R. and D. Hilbert: \emph{Methoden der mathematischen Physik
II,} Springer, Berlin-Heidelberg-New York 1968.

3.\qquad Dzyadyk, V. K.: \emph{Introduction to the Theory of Uniform
Approximation of Functions by Polynomials,} Nauka, Moscow 1977.

4.\qquad Helms, L. L.: \emph{Einfiihrung in die Potentialtheorie, de Gruyter},
Berlin-New York 1973.

5.\qquad Johnen, H., and K. Scherer: On the equivalence of the A"-functional
and moduli of continuity and some applications, in: \emph{Constr. Theory of
Funct. of Several Variables, }Lecture Notes Math. 571, pp. 119-140, Springer,
Berlin-Heidelberg-New York (1977).

6.\qquad Kounchev, O. I.: Harmonicity moduli and Jackson type theorems for the
approximation through polyharmonic functions (preprint).

7.\qquad Meinardus, G.: \emph{Approximation of Functions: Theory and Numerical
Methods,} Springer, Berlin-Heidelberg-New York 1967.

8.\qquad Nicolescu, M.: \emph{Opera Matematica. Functii Poliarmonice,} Editura
Academiei, Bucuresti 1980.

9.\qquad Nikolskii, S. M.: \emph{Approximation of Functions of Several
Variables and Imbedding Theorems, }Springer, Berlin-Heidelberg-New York 1975.

\bigskip

Received September 16, 1991.

\end{document}